\documentclass[12pt,reqno]{amsart}
\usepackage{verbatim}
\usepackage{amscd, amsfonts}
\usepackage[mathscr]{eucal}

\begin{document}

\numberwithin{equation}{section}
\renewcommand{\theequation}{\thesection.\arabic{equation}}
\setcounter{secnumdepth}{2}

%
\newcommand{\Ac}{{\mathcal A}}
\newcommand{\Bc}{{\mathcal B}}
\newcommand{\bt}{\tilde{b}}
\newcommand{\bfF}{\bf{F}}
\newcommand{\bfG}{\bf{G}}
\newcommand{\BH}{\Bc_H}

\newcommand{\Dc}{{\mathcal D}}
\newcommand{\Fc}{{\mathcal F}}
\newcommand{\Gma}{{\gamma}}
\newcommand{\Hc}{\mathscr{H}}
\newcommand{\htil}{\tilde{h}}

\newcommand{\Lap}{\Delta}
\newcommand{\Mcl}{{\mathcal M}}

\newcommand{\nab}{\nabla}
\newcommand{\Om}{\Omega}
\newcommand{\Omb}{\overline{\Omega}}

\newcommand{\pal}{\partial}
\newcommand{\Pc}{\mathbb{P}}
\newcommand{\psit}{\tilde{\psi}}

\newcommand{\Rb}{\overline{\R}}
\newcommand{\R}{\mbox{$\mathbb R$}}

\newcommand{\Sc}{{\mathcal S}}
\newcommand{\sg}{\sigma}
\newcommand{\sgt}{\tilde{\sigma}}
\newcommand{\vap}{\varphi}

\newcommand{\Wc}{{\mathcal W}}

\newcommand{\barr}{\begin{eqnarray}}
\newcommand{\bc}{\begin{center}}
\newcommand{\beq}{\begin{equation}}
\newcommand{\bpf}{\begin{proof} \quad}
\newcommand{\btm}{\begin{thm}}

\newcommand{\earr}{\end{eqnarray}}
\newcommand{\ec}{\end{center}}
\newcommand{\eeq}{\end{equation}}
\newcommand{\epf}{\end{proof}}
\newcommand{\etm}{\end{thm}}

\newcommand{\ang}[1]{\langle#1\rangle}
\newcommand{\bdy}{\partial \Omega}
\newcommand{\Cinty}{C^{\infty} }

\newcommand{\deq}{:= }
\newcommand{\deqs}{\ :=\ }

\newcommand{\eqs}{\ =\ }
\newcommand{\geqs}{\ \geq \ }
\newcommand{\leqs}{\ \leq \ }
\newcommand{\mns}{\, - \, }
\newcommand{\pls}{\, + \, }
\newcommand{\plms}{\ + \ }


\newcommand{\foral}{\qquad \mbox{for all} \quad }
\newcommand{\foreach}{\qquad \mbox{for each} \quad }

\newcommand{\wrt}{ \mbox{with respect to}  }
\newcommand{\xand}{\quad \mbox{and} \quad }
\newcommand{\xfor}{ \quad \mbox{for} \ }
\newcommand{\xiff}{\ \mbox{if and only if} \ }

\newcommand{\xon}{\qquad \mbox{on} \quad}
\newcommand{\xor}{\qquad \mbox{or} \quad}
\newcommand{\xthen}{\quad \mbox{then} \quad}
\newcommand{\xwhen}{\qquad \mbox{when} \quad}
\newcommand{\xwith}{\qquad \mbox{with} \quad}


\newcommand{\Cinfty}{C^{\infty}}

\newcommand{\delj}{\delta_j}
\newcommand{\Dbj}{\Dnu b_j}
\newcommand{\Dej}{ \Dnu e_j    }
\newcommand{\DFu}{D \Fc(u)}

\newcommand{\Div}{\mathop{\rm div}\nolimits}
\newcommand{\divF}{\Div \bfF}
\newcommand{\divG}{\Div \bfG}

\newcommand{\Dnu}{D_{\nu}}
\newcommand{\dsg}{\, d \sigma}
\newcommand{\dsgt}{\; d \sgt}

\newcommand{\ej}{e_j}
\newcommand{\ejh}{\hat{\ej}}
\newcommand{\fjh}{\hat{f}_j}
\newcommand{\fkh}{\hat{f}_k}

\newcommand{\gamH}{\gamma_H}
\newcommand{\gjh}{\hat{g}_j}
\newcommand{\hjh}{\hat{h}_j}
\newcommand{\HcM}{\Hc_M}
\newcommand{\kjh}{\hat{k}_j}

\newcommand{\Fnu}{\bfF \cdot \nu}
\newcommand{\Gjhat}{\hat{G}_j}

\newcommand{\gradchi}{\nabla \chi}
\newcommand{\gradphi}{\nabla \varphi}
\newcommand{\gradpsi}{\nabla \psi}
\newcommand{\gradu}{\nabla u}
\newcommand{\gradv}{\nabla v}

\newcommand{\Iby}{\int_{\bdy} \; }
\newcommand{\IOm}{\int_{\Om} \; }

\newcommand{\Keta}{K_{\eta}}
\newcommand{\lamo}{\lambda_1}
\newcommand{\Lapbj}{\Delta \, b_j}
\newcommand{\Lapchi}{\Delta \chi}
\newcommand{\Lapsi}{\Delta \psi}

\newcommand{\Lapu}{\Delta \, u}
\newcommand{\Lapv}{\Delta \, v}

\newcommand{\mut}{\tilde{\mu}}
\newcommand{\n}[1]{\left\vert#1\right\vert}
\newcommand{\nm}[1]{\left\Vert#1\right\Vert}

\newcommand{\PcH}{{\Pc}_H}
\newcommand{\Omt}{\Om \times \Om}
\newcommand{\opal}{\oplus_{\pal}}
\newcommand{\qj}{q_j}
\newcommand{\qone}{q_1}

\newcommand{\Rn}{{\R}^N}
\newcommand{\RN}{{\R}^N}
\newcommand{\ra}{\rightarrow}
\newcommand{\ROm}{R_{\Om}}
\newcommand{\ROmM}{R_{\Om M}}

\newcommand{\sumM}{\sum_{j=1}^M}
\newcommand{\sumo}{\sum_{j=1}^{\infty}}
\newcommand{\sumz}{\sum_{j=0}^{\infty}}

\newcommand{\ujh}{\hat{u}_j}
\newcommand{\ut}{\tilde{u}}
\newcommand{\uzh}{\hat{u}_0}

\newcommand{\vjh}{\hat{v}_j}
\newcommand{\vt}{\tilde{v}}
\newcommand{\vzh}{\hat{v}_0}


\newcommand{\HLap}{H(\Delta, \Om)}
\newcommand{\Hone}{H^1(\Om)}
\newcommand{\Hmone}{H^{-1}(\Om)}
\newcommand{\Hstar}{\Hone^*}
\newcommand{\Honet}{\Hone \times \Hone}
\newcommand{\Hzone}{H_{0}^1(\Om)}
\newcommand{\HzLap}{H_0(\Delta,\Om)}
\newcommand{\HzzLap}{H_{00}(\Delta,\Om)}
\newcommand{\Lp}{L^p (\Om)}
\newcommand{\Lq}{L^q (\Om)}
\newcommand{\Lt}{L^2 (\Om)}


\newcommand{\BHarm}{\Bc \Hc(\Om)}
\newcommand{\Btom}{L^2_H(\Om)}
\newcommand{\Cc}{C_c^1 (\Om)}

\newcommand{\Harm}{\Hc (\Om)}
\newcommand{\Hharm}{{\Hc}^{1/2}(\Om)}
\newcommand{\Hsharm}{\Hc^s(\Om)}
\newcommand{\Hspharm}{{\Hc}^{s+1/2}(\Om)}


\newcommand{\Hdiv}{H(\Div, \Om)}
\newcommand{\LtN}{L^2 (\Om; \RN)}

\newcommand{\Hhby}{H^{1/2}(\bdy)}
\newcommand{\Hmhby}{H^{-1/2}(\bdy)}

\newcommand{\Hmsby}{H^{-s}(\bdy)}
\newcommand{\Hsby}{H^s(\bdy)}

\newcommand{\Lpb}{L^p (\bdy, \dsg)}

\newcommand{\Ltby}{L^2 (\bdy, \dsg)}
\newcommand{\Lttby}{L^2 (\bdy, \dsgt)}

%
\newtheorem{thm}{Theorem}[section]
\newtheorem{cor}[thm]{Corollary}
\newtheorem{cond}{Condition}
\newtheorem{lem}[thm]{Lemma}
\newtheorem{prop}[thm]{Proposition}

\title[  S.V.D. of the Poisson kernel]
{The S.V.D. of the Poisson Kernel}
\author[Auchmuty]{Giles Auchmuty}

\address{Department of Mathematics, University of Houston, 
Houston, Tx 77204-3008 USA}
\email{auchmuty@uh.edu}

\thanks{The author gratefully acknowledges research support by NSF award DMS 11008754. 
 }

\date{April 19, 2016}

\begin{abstract} 
This paper describes a singular value decomposition (SVD) for the Poisson kernel associated with the Laplacian on 
bounded regions  $\Om$ in $\RN, \ N \ge 2$. 
The singular functions and singular values are related to certain  Steklov eigenvalues and eigenfunctions of the biharmonic 
operator on $\Om$.
This eigenproblem and its properties are studied on an appropriate space. 
This enables a description of the Bergman harmonic projection and an  orthonormal basis of the real  harmonic Bergman space 
$\Btom$ is found.
 A reproducing kernel  for $\Btom$ is constructed and also an  orthonormal basis of the space $\Ltby$.
The Poisson kernel may be regarded as   the harmonic extension operator from $\Ltby$ to $\Btom$ and  has an explicit 
spectral representation that yields the SVD of the Poisson kernel. 
The singular values of the Poisson kernel are related to the eigenvalues of the DBS eigenproblem. 
This enables the description of optimal finite rank approximations of the Poisson kernel with error estimates. 
Explicit spectral formulae for the normal derivatives of eigenfunctions for the Dirichlet Laplacian on $\bdy$ are found and
 used to identify a constant in an inequality of Hassell and Tao. 
\end{abstract}
\maketitle

\section{Introduction}\label{s1}

This paper describes some new  representation results for    harmonic functions on a bounded region 
$\Om$ in $\RN; \, N \geq 2. $
In particular an explicit description of the Reproducing Kernel for the harmonic Bergman space $\Btom$ and 
the SVD of the Poisson kernel will be obtained.
In a earlier paper \cite{Au3} the author  used harmonic Steklov eigenfunctions to represent  reproducing kernels on the family $\Hsharm$ of harmonic functions.
These spaces are modeled on Sobolev spaces and the Steklov eigenfunctions were not, in general,   $L^2$ orthogonal - though their gradients are. 

Here we shall show how similar methods may be used to construct an  orthonormal basis of  the harmonic Bergman space $\Btom$ for a general class  of bounded regions.
This basis is obtained by constructions involving the eigenfunctions of the {\it Dirichlet Biharmonic Steklov  (DBS) eigenproblem}.
G. Fichera \cite{Fic} in 1955 showed  that the norm of the Poisson kernel as a map of $\Ltby$ to $\Lt$ is a function of the
smallest eigenvalue of this  problem. 
Here an expression for this harmonic extension operator $E_H$ in terms of these DBS eigenfunctions and eigenvalues is found.
This eigenfunction expansion provides  an SVD for the Poisson kernel that holds under  mild boundary  regularity requirements.

The results obtained here are found under weaker hypotheses on the boundary $\bdy$ of the region than have previously 
been used for studies of these problems. 
This is possible due to the use of special Sobolev spaces related to the Laplacian that are introduced in section \ref{s3}.
A crucial result is theorem \ref{T3.1} that proves a continuity result for  the normal derivative operator $\Dnu$.

Results about some Steklov eigenproblems related to the Laplacian are described in sections \ref{s4} and \ref{s5}.
A nice summary of recent results on  the DBS eigenproblem may be found  in chapter 3 of Gazzola-Grunau-Sweers \cite{GGS}. 
In particular properties of the spectrum were obtained by Ferrero, Gazzola and Weth \cite{FGW} and properties of the first 
eigenvalue were studied in Bucur, Ferrero and Gazzola \cite{BFG}. 
Here the DBS eigenproblem is studied  in the space $\HzLap$ which may be different to the space used in previous treatments. 
The results follow from an algorithm to construct an explicit sequence   of Dirichlet Biharmonic Steklov (DBS) eigenfunctions 
that yield an orthonormal  basis $\BH$ of $\Btom$.
 This construction is described  in section \ref{s5}.
 
 In section \ref{s6},  the {\it Bergman harmonic projection} $\PcH$ of $\Lt$ onto $\Btom$ is first defined. 
 Some differences between the  harmonic projection on $\Hone$  used in elliptic PDEs  and the Bergmann harmonic
 projection  are  described.
 Then the sequence of DBS eigenfunctions is used to  construct  an orthonormal  basis of $\Btom$.
 A representation result   for the Bergman harmonic projection $\PcH$  of $\Lt$ onto $\Btom$ in terms of the  DBS 
 eigenfunctions is proved as theorem \ref{T6.1}.

A formula for the {\it  reproducing kernel} (RK)  for $\Btom$ follows as  corollary \ref{C6.3}. 
This RK may be viewed as a Delta function  on the class of harmonic functions as \eqref{eDel} holds. 
The reproducing kernel described here appears to be  a spectral representation of a reproducing kernel constructed  
by J.L. Lions in  \cite{L1}  using  control theory methods.
Lions showed that there is a reproducing kernel for $\Btom$ that is  a perturbation of the fundamental solution 
of the biharmonic operator;  the perturbation depending on the region $\Om$. 
Subsequently, Englis et al in \cite{ELPP} and J.L. Lions in \cite{L2}, have studied the construction of other 
Reproducing Kernels  for various classes of harmonic and other elliptic operators on bounded regions.  
 
The SVD of the classical Poisson integral  operator regarded as a linear transformation  from $\Ltby$ to the harmonic
Bergman space $\Btom$  is described in section \ref{s7}. 
A spectral representation of the harmonic extension operator is first described and shown to be compact.
When this operator is regarded as an integral operator then its kernel is the Poisson kernel and we observe that 
the representation  may be regarded as  a SVD.
The singular vectors are the orthonormal bases $\BH$ and $\Wc$ involving  DBS eigenvalues and eigenfunctions and
the singular values are related to the DBS eigenvalues.
Moreover  associated  finite rank approximations of the Poisson operator  have error estimates  depending on 
appropriate  DMS eigenvalues.
See theorem \ref{T7.2} and the error estimates for finite rank approximations of the Poisson  operator in 
theorem \ref{T7.3}.

In section \ref{s8}  an explicit formula for the normal derivative of Dirichlet Laplacian eigenfunctions on $\Om$ is found.
This provides a quantification of a constant described by Hassell and Tao \cite{HT} for the 2-norms of such 
eigenfunctions in the case where the domain is a nice bounded region of $\RN$.

The results here are stated under  a weak  regularity condition (B2) on the boundary $\bdy$.  
This condition has been the subject of recent interest as it is related to phenomena that arise in the study of 
biharmonic boundary  value problems.
Some comments about these issues may be found in  section 2.7 of the monograph of Gazzola, Grunau and Sweers \cite{GGS} 
 including a description  of some apparent "paradoxes".

\vspace{1em}

\section{Definitions and Notation.}\label{s2}

   A region  is a non-empty, connected, open subset of $\Rn$. Its closure  is denoted  $\Omb$ and its boundary is   
   $\bdy \deq \Omb\setminus \Om$. 
   A standard assumption about the region is the following. 

\noindent{\bf (B1):}  {\it $\Om$ is a bounded region in  $\Rn$ and its boundary  $\bdy$  is the union of a finite number of 
disjoint closed Lipschitz  surfaces; each surface having finite surface area.}

When this holds there is an outward unit normal $\nu$ defined at $\sigma \; a.e.$  point of $\bdy$. 
The  definitions and terminology of Evans and Gariepy \cite{EG}   will be followed except  that $\sigma, \dsg$,  respectively, 
will represent Hausdorff   $(N-1)-$dimensional measure and integration with respect to this measure. 
All functions in this paper will take values in  $\Rb \deq [-\infty,\infty]$ and derivatives should be taken in a weak sense.

 The real Lebesgue spaces $\Lp$ and $\Lpb, \ 1 \le p \le \infty$ are  defined in the standard manner and have the usual norms
  denoted by  ${\nm{u}_p}$ and  $ {\nm{u}}_{p,\bdy}$.
   When $p = 2$, these spaces will be Hilbert  spaces with  inner products  
\[ \ang{u , v}  \deqs \int_{\Om} \  u(x) \; v(x) \ dx \qquad \mbox{and} \quad
 \ang{u , v}_{\bdy}  \deqs |\bdy|^{-1} \ \int_{\bdy} \  u \; v \ \dsg.  \]
 
 Let $\Hone$ be the usual real Sobolev space of  functions on $\Om$.  
 It is a real Hilbert space under  the standard $H^1-$ inner product
\begin{equation}\label{ip1}
{[u , v]}_1  \deq \IOm \ [ u(x)\, v(x) \  + \  \gradu(x)  \cdot  \gradv(x)] \  dx.
\end{equation}
Here  $\gradu$ is the gradient of the function $u$ and the associated norm   is denoted   $ {\nm{u}}_{1,2}$.

The region $\Om$ is said to satisfy {\it Rellich's theorem} provided  the imbedding  of  $\Hone$ into $\Lp$ is compact for 
$1 \leq p < p_S$ where $p_S := \, 2N/(N-2); \ N \geq 3$, or $p_S = \infty$ when $N = 2$. 
There are a number of different criteria on $\Om$ and $\bdy$ that imply this result. 
When (B1) holds it is theorem 1 in section 4.6 of \cite{EG}; see also Amick  \cite{Am}.  
DiBenedetto \cite{DB}, in theorem 14.1 of chapter 9 shows that the result holds when $\Om$ is bounded and satisfies a "cone property". 
Adams and Fournier give a  thorough  treatment of conditions for this result in chapter 6 of \cite{AF} and show that it also
holds for some classes of unbounded regions.

When (B1) holds, then the trace of a Lipschitz continuous function on  $\Omb$ to $\bdy$ is continuous and there is a continuous  extension of this map to $W^{1,1}(\Om)$.
This linear map $\Gma$ is called the  trace on  $\bdy$ and  each $\Gma (u)$  is  Lebesgue integrable with respect to  $\sigma$; see \cite{EG}, section 4.2 for details. 
In particular, when  $\Om$ satisfies (B1), then the {\it Gauss-Green} theorem  holds in the form
\beq \label{GG} 
\IOm   u(x) D_j \, v(x) \  dx \eqs  \Iby u \, v \, \nu_j  \dsg \mns  \IOm   v(x) D_j \, u(x) \  dx \quad \mbox{for} \ 1 \le j \le N.
\eeq
and all $u, v$ in $\Hone$. Often, as here, $\Gma$ is omitted in boundary integration.

The region $\Om$ is said to satisfy a {\it compact trace theorem} provided  the  trace mapping $\Gma : \Hone \ra \Ltby$ is compact. 
  Evans and Gariepy \cite{EG}, section 4.3 show that $\gamma$ is  continuous when $\bdy$ satisfies (B1). Theorem 1.5.1.10 of Grisvard  \cite{Gri} proves an inequality that implies the compact trace theorem when $\bdy$  satisfies (B1). This inequality is also proved in \cite{DB}, chapter 9, section 18 under stronger regularity   conditions on the boundary. 

We will generally  use  the following  equivalent inner product on $\Hone$
\begin{equation}\label{ip2}
{[u , v]}_{\pal}  \deq \IOm \  \gradu  \cdot \nabla v \  dx \pls  \Iby \  u \; v \ \dsgt.
\end{equation}
 The related  norm is denoted $ {\nm{u}}_{\pal}$. $\sgt$ is the normalized  surface area measure  defined by  
 $\sgt(E) \deqs |\bdy|^{-1} \sg(E)$ where $|bdy| := \sg(\bdy)$ is the surface measure of the boundary. 
  The proof that this  norm is equivalent to the usual $(1,2)-$norm on $\Hone$ when  (B1) holds is Corollary 6.2 of \cite{Au1} 
  and also is part of theorem 21A of \cite{Z1}. 

When $\bfF \in {\rm \LtN}$ and there is a function $\vap \in \Lt$ satisfying
\beq \label{e2.11}
\IOm u \, \vap \ dx \eqs \IOm \gradu \cdot \bfF \ {\rm dx} \foral {\it u \in \Cc }\eeq
then we say that $\divF \deqs \vap$ is the divergence of $\bfF$.
The class of all $L^2-$vector fields on $\Om$ whose divergence is in $\Lt$ is denoted $\Hdiv$ and is a real Hilbert space with
the inner product
\beq \label{e2.12}
[\bfF, \bfG]_{div} \deqs \IOm [ \bfF \cdot \bfG \pls \divF \, \divG ] \ {\rm dx.} \eeq

The results described here depend on techniques and results of variational calculus. 
Relevant notations and definitions are those of  Attouch, Buttazzo and Michaille \cite{ABM}. 

\vspace{1em}


\section{The spaces $\HLap$ and  $\HzLap$. } \label{s3}

Henceforth  the region $\Om$  is assumed to satisfy (B1).
Define $\HLap$ to be the subspace  of all functions $u \in \Hone$ with $\gradu \in \Hdiv$. 
Write $\Delta u \deqs \Div (\gradu)$ so $\Delta$ is   the usual Laplacian. 
$\HLap$ is a real Hilbert space $\wrt$ the inner product
\beq \label{e3.1}
[u, v]_{\pal,\Delta} \deqs [u,v]_{\pal} \pls \IOm \, \Lapu \, \Lapv  \ dx   . \eeq 

A function $u \in \Lt$ is said to be {\it harmonic} on $\Om$ provided
\beq \label{e3.2}
\IOm \ u \ \Delta v \ dx \eqs 0 \foral v \in C_c^2(\Om) . \eeq
Thus a function  $u \in \Hone$ will be harmonic on $\Om$ provided
\beq \label{e3.3}
\IOm \ \gradu \cdot \gradv \ dx \eqs 0 \foral v \in \Hzone. \eeq
Let $\Harm$ be the class of all $H^1-$ harmonic functions on $\Om$.
The following result has been used in a variety of ways in some preceding papers, 
\cite{Au1} and \cite{Au2}, that study  other issues.
Here a different statement and a direct proof is provided for   completeness.

\begin{lem} \label{L3.1}
Suppose that $\Om$ satisfies (B1). Then there are closed subspace $\Hzone, \Harm$ of
$\Hone$ and projections $P_0, P_H$ onto these spaces  such that 
\beq \label{e3.5}
u \eqs P_0 u \pls P_H u \foral u \in \Hone. \eeq
Moreover $\gamma(u) = \gamma(P_H u)$ and  $[P_0 u , P_H u]_{\pal} \eqs 0$ for all $u \in \Hone$. 
\end{lem}
\bpf
Given $u \in \Hone$, consider the variational problem of minimizing
\[ \Fc(v) \deqs \| \, v \mns u \, \|_{\pal}^2 \quad \mbox{over} \quad v \in \Hzone. \]
This problem has a unique minimizer $u_0 \in \Hzone$ as $\Fc$ is convex, coercive and 
continuous on $\Hzone$. 
Evaluation of  the G-derivative of $\Fc$ implies that   the minimizer satisfies 
\[ D \Fc (u_0) (v) \eqs 2 \IOm \ \nabla(u_0 - u) \cdot \gradv \ dx \eqs 0 \foral v \in  \Hzone. \]

That is $u_h := u - u_0$ is harmonic on $\Om$. 
Define $P_0 u = u_0$ and $P_H u = u_h$, these are continuous maps into $\Hzone, \Harm$
respectively. 
These are projections with closed  range  from corollary 3.3 of Auchmuty \cite{Au}.
Since $\gamma(u_0) = 0$ one has $\gamma(u) = \gamma(u_h)$ and the orthogonality 
follows from the extremality condition above.    \epf

This lemma provides a $\pal-$orthogonal decomposition of $\Hone$ and the operator $P_H$ 
defined here is the standard harmonic projection of $H^1$ functions.

Define $\HzLap$ to be the range of $P_0$ when restricted to $\HLap$.
It is a closed subspace of $\HLap$ and the orthogonal decomposition 
\beq \label{e3.7}
\HLap \eqs \HzLap \oplus_{\pal, \Delta} \Harm  \eeq
holds with respect to the inner product \eqref{e3.1}. 

The following theorem shows that when $u \in \HzLap$, the boundary flux $\Dnu u$ has further regularity. 

\btm \label{T3.1}
 \quad Suppose that (B1) holds and  $u \in \HzLap$. 
 Then $\Dnu u$ is in $\Ltby$  and there is a $C_{\Om}$ such that 
 $\| \Dnu u\|_2 \leqs C_{\Om} \ \| \Delta u \|_2 $ for all $u \in \HzLap$.      \etm
 \vspace{-1em}
\bpf  When (B1) holds let $v \eqs v_0 + v_h$ with $v_0 \in \Hzone$ and $v_h  \in \Harm$ be a decomposition
of $v \in \HLap$ as in lemma \ref{L3.1}.
When $u \in \HzLap$, Green's formula for Sobolev functions on $\Om$ becomes
\beq \label{e3.8}
\IOm [v_h \, \Delta u - u \, \Delta v_h] \ dx \eqs \Iby \, \gamma(v_h) \, \Dnu u \, \dsg. \eeq

Since   $v_h$ is harmonic on $\Om$ and $\gamma(v_h) = \gamma(v)$, this becomes
\[ |\bdy| \, |\ang{\gamma(v), \Dnu u}_{\bdy}| \leqs \| v_h \|_{2,\Om}  \  \| \Delta u \|_2  \foral v \in \Hone. \]

From theorem 6.3 of \cite{Au3},  $\Hc^{1/2}(\Om)$ and $\Ltby$ are isometrically isomorphic, so 
\[  \| \Dnu u \|_{2, \bdy}   \eqs \sup_{\|v\|_{2,\bdy} \leq 1} |\ang{v, \Dnu u}_{\bdy}| \leqs  
\sup_{\|v_h\|_{1/2,\Om} \leq 1} \|v_h\|_{2,\Om} \,   \| \Delta u \|_2   \]
\[ \leqs C_{\Om} \ \| \Delta u \|_2   \xwith   C_{\Om} \deqs \sup_{\|v_h \|_{1/2,\Om} \leq 1}    \| v_h \|_{2,\Om}. \]
This $C_{\Om}$  is finite and attained as the imbedding of $\Hc^{1/2}(\Om)$ into $\Lt$ is compact. 
\epf

Now consider  the inner product on $\HzLap$ defined by
\beq \label{e3.15}
[u, v]_{\Delta} \deqs  \IOm \, \Lapu \, \Lapv  \ dx   . \eeq 
The following inequality shows that this generates an equivalent norm to that of  $(\pal,\Delta)$.

\begin{lem} \label{L3.3}
 \quad Suppose that (B1) holds, $u \in \HzLap$ and $\lamo$ is the first  eigenvalue of the Dirichlet 
 Laplacian on $\Om$.   Then
\beq \label{e3.17}
\| u \|_{\Delta}^2 \leqs  \| u \|_{\pal,\Delta}^2 \leqs (1\pls \frac{1}{\lamo}) \ \| u \|_{\Delta}^2  \foral  u \in \HzLap. \eeq
\end{lem}
\bpf 
The first inequality is trivial, while the second follows from the fact that 
\[ \IOm |\Lapu|^2 dx \geqs \lamo \IOm |\gradu|^2 dx \foral  u \in \HzLap.   \]   
\epf

When $\bdy$ satisfies further regularity conditions, it is well known that $\HzLap \eqs H^2(\Om) \cap \Hzone$.
This is proved in Evans \cite{Ev} chapter 6, section 6.3.2 when $\bdy$ is $C^2$. 
Adolfsson \cite{Ado} has shown that this holds when  $\bdy$ is bounded, Lipschitz and  satisfies a uniform outer ball 
condition.

For this paper a slightly stronger assumption than (B1) about the region $\Om$ is needed namely;

\noindent{\bf (B2):} \quad   {\it $\Om$ is a bounded region with a boundary $\bdy$ for which  (B1) holds and 
$\Dnu$ is a compact mapping of $\HzLap$ into $\Ltby$. }

This condition has been verified under various regularity conditions on the boundary $\bdy$.
Necas \cite{Nec}  chapter 2, theorem 6.2 has shown that (B2) holds when $\Om$ is Lipschitz and satisfies a uniform outer ball condition.
Grisvard \cite{Gri} chapter 1.5 has a further discussion of this. 
 (B2) also holds when each component of the boundary $\bdy$ is a $C^2-$manifold.
More literature about this is described  in section 2.7 of \cite{GGS} where  some related "paradoxes" are  presented. 


\vspace{1em}
\section{Harmonic Steklov Representations and Boundary Traces.  } \label{s4}

The  methods used here depend on  results about  boundary traces described in some earlier papers of the author.  
In particular, the spectral characterization of trace spaces described in Auchmuty \cite{Au2} and results about spaces of harmonic functions proved in \cite{Au3} will be used. 
For convenience some of these results are summarized below. 
Henceforth $\Om$ is a bounded region in $\RN$ satisfying (B2). 

A function $s \in \Hone$ is said to be a {\it harmonic Steklov eigenfunction} provided it is a non-zero solution of
\beq \label{e4.05}
\IOm \nabla s \cdot \gradv \, dx \eqs \delta \ \Iby s \, v \dsgt \foral v \in \Hone. \eeq
When this holds then $\delta$ is the associated Steklov eigenvalue. 

Let $\Sc \deqs \{s_j : j \geq 0 \}$ be a maximal orthogonal sequence of harmonic Steklov eigenfunctions as described 
in \cite{Au2}. 
Assume that they are normalized so that their boundary traces are $L^2-$orthonormal; $\ang{s_j,s_k}_{\bdy} = \delta_{jk}$
for all $j,k$.
A function $f \in \Lttby$ has the usual representation 
\beq \label{e4.1}
f(x) \eqs  \sum_{j=0}^{\infty} \  \fjh s_j(x) \xon \bdy  \xwith \fjh := \ang{f,s_j}_{\bdy} \eeq
with respect to this basis. 
Here $\fjh$ is called the j-th Steklov coefficient of $f$ and \eqref{e3.1} is  called the  Steklov representation of $f$.

When $f \in \Lttby$ then the function   $Ef \in \Lt$ defined by 
\beq \label{e4.2}
Ef(x)  \deqs \sumz \fjh \, s_j(x) \eeq
is a harmonic function on $\Om$.

This $f$ is defined to be  in the trace space $\Hsby$ with $s  \geq 0$ provided its Steklov coefficients   satisfy
\beq \label{e4.3}
\| f \|_{s,\bdy}^2  \deqs \sum_{j=0}^{\infty} \ (1 \pls \delj)^{2s} \  |\fjh|^2  \  < \ \infty . \eeq  
Thus $s=0$ is the usual Lebesgue space $\Ltby$ from Parseval's identity.
When $f \in \Hsby$, then $Ef$ is in a space $\Hspharm$ and moreover $E$ is an isometric isomorphism of these spaces.
See theorems 6.2 and 6.3 in \cite{Au3} for full statements and proofs of this. 

A linear functional $G$ is in the dual space $\Hmsby$ provided  G has Steklov coefficients $\Gjhat$ and there is a constant 
C such that  
\beq \label{e4.4}
G(f) \deqs \ang{G,f}_{\bdy} \deqs \sum_{j=0}^{\infty} \  \Gjhat \fjh  \leqs C \ \|f \|_{s,\bdy}  \foral  f \in \Hsby.
\eeq
This pairing  of $\Hsby$ and $\Hmsby$ extends the usual $L^2-$ inner product.
 That is, functionals in  $\Hmsby$ have well-defined Steklov representations so they may be regarded as  {\it generalized functions} 
 on $\bdy$ .
 
When $u \in \Hone$, then its boundary trace $\gamma(u)$ will be in $\Hhby$ and its normal derivative $\Dnu u$ is a generalized
function in $\Hmhby$. Suppose
\beq \label{e4.5}
\gamma(u)(x)  \eqs \sum_{j=0}^{\infty}  \ujh \, s_j(x),  \xthen\ \Dnu u \eqs   
\sum_{j=1}^{\infty}  \ \delj \, \ujh \, s_j. \eeq 

Note that the inner product on $\Hhby$ associated with \eqref{e4.3} is
\beq \label{e4.7}
[u,v]_{1/2,\bdy} \deqs \uzh \, \vzh \pls \ang{\Dnu u, v}_{\bdy} \eeq
and this expression is symmetric in $u,v$.
When $u \in \HzLap$, then theorem \ref{T3.1} implies the last term here is a standard boundary integral as $\Dnu u \in \Ltby$.

When (B1) holds and $u, v \in \HLap$, then a classical  Green's identity becomes
\beq \label{e4.9}
 \IOm [ \, u \, \Lapv \mns v\,  \Lapu \, ] \ dx \eqs   |\bdy| \left[ \, \ang{\Dnu v, u}_{\bdy} \mns  \ang{\Dnu u, v}_{\bdy} \, \right]
 \eeq
where the terms on the right hand side are defined by  pairings  of the form
\beq \label{e4.11}
\ang{\Dnu u, v}_{\bdy} \deqs \sum_{j=1}^{\infty}  \ \delj \ \ujh \, \vjh   \eeq
that extend standard boundary integrals. As a consequence one sees that
\beq \label{e4.12}
\IOm u \, \Delta v \, dx \eqs \IOm v \, \Delta u \, dx \foral u, v \in \HzLap. \eeq

\vspace{1em}
\section{Dirichlet  Biharmonic Steklov Eigenproblems on $\Om$.} \label{s5}

In this section we shall describe some properties of solutions of the {\it Dirichlet Biharmonic Steklov} (DBS) eigenproblem
on a region $\Om \subset \RN$ that satisfies  (B2).
This problem is described in section 5.1 of Kuttler and Sigillito \cite{KS} and more recently in section 3.3 of 
Gazzola, Grunau and Sweers   \cite{GGS}. 
Note that (B2)  is  weaker  than the requirements on the domain in   the analyses of \cite{GGS} and others. 

Our particular aim is to construct  an orthonormal basis of $\HzLap$ using the  framework of Auchmuty \cite{Au4}. 
This involves the solution of a sequence of constrained  variational principles - which are not your standard
principles involving Rayleigh quotients and equality constraints.

A function $b \in \HLap$ is said to be (weakly) {\it biharmonic} provided
\beq \label{e5.3}
\IOm \Delta b \, \Delta v \ dx \eqs 0 \foral v \in C_c^2(\Om). \eeq 
The {\it DBS eigenproblem}  is to find nontrivial solutions $(q, b) \in \R \times \HzLap$ of the system
\beq \label{e5.1}
\IOm \Delta b \, \Delta v \ dx \eqs q \ \Iby \, \Dnu b \, \Dnu v \dsg    \foral v \in \HzLap. \eeq

This is a weak version of the problem of finding biharmonic functions $b \in \HzLap$ that satisfy the 
boundary conditions
\beq \label{e5.2}
b \eqs \Delta b \mns q \, \Dnu b \eqs 0 \xon \bdy. \eeq
Here $q$ is the DBS eigenvalue and this is a Steklov eigenprobem as the eigenvalue appears only in the
boundary condition.

Take $V = \HzLap$, then \eqref{e5.1} has the form of the problem studied in Auchmuty \cite{Au4} with the notation,
\beq   \label{e5.4}
a(u,v) \deqs \IOm \Delta u \, \Delta v \ dx, \quad m(u,v) \deqs \Iby \, \Dnu u\, \Dnu v \dsg  \xand \lambda := q. 
\eeq

Functions $u,v$  in $\HzLap$ are said to be $\Delta-$orthogonal (resp. m-orthogonal) provided $a(u,v) = 0, 
(m(u,v) = 0)$. 
When $b_1, b_2$ are two DBS eigenfunctions satisfying \eqref{e5.1}, then they will be $\Delta-$orthogonal if and only
 if they are  m-orthogonal. 
First note that if \eqref{e5.1} has a nontrivial solution $(q,b)$ then by letting $v=b$ it follows that  $ q \geq 0$.
If $q = 0$ then $\Delta b \equiv 0$ on $\Om$ so by uniqueness $b \equiv 0$.
Thus all DBS eigenvalues must be  strictly positive.

To find the smallest DBS eigenvalue, let $C_1$ be the closed unit ball in $\HzLap$ and 
consider the problem of maximizing
\beq \label{Max1}
\Mcl(u) \deqs \Iby \ |\Dnu u|^2 \dsg   \quad \mbox{subject to} \qquad \| u \|_{\Delta}^2 \leqs 1. \eeq
Define  $\beta_1 \deqs \sup_{u \in C_1} \ \Mcl(u), $ then the following result generalizes  the existence results of 
theorem 3.17 of \cite{GGS} and \cite{BFG}. 
It requires weaker assumptions on the boundary $\bdy$ and the solutions are in a different space.

\btm \label{T5.1}
Assume that (B2) holds, then there are functions $\pm b_1 \in C_1$ that maximize $\Mcl$ on $C_1$.
These functions are non-trivial solutions of \eqref{e5.1} associated with the smallest positive eigenvalue 
$q_1 = 1/{\beta_1}$ of the DBS eigenproblem and
\beq \label{e5.5}
\IOm \ |\Delta u|^2 \ dx \geqs q_1 \ \Iby  |\Dnu u|^2 \dsg  \foral u \in \HzLap. \eeq \etm
\bpf
$C_1$ is weakly compact in, and $\Mcl$ is weakly continuous on  $\HzLap$ so $\Mcl$ attains its supremum on $C_1$. 
If $b_1$ is such a maximizer so also is $-b_1$ as $\Mcl$ is even. 
Let $I_1(u)$ be the indicator functional of $C_1$, then from part (ii) of  Theorem 9.5.5 of Attouch, Buttazzo and Michaille \cite{ABM}
the maximizers are solutions of the inclusion $ \, 0 \ \in \ D \Mcl(b) \pls \pal I_1 (b)$. 
Here $D$ is a G-derivative and $\pal$ denotes the subdifferential.
In proposition 9.6.1, it is shown that $\pal I_1(u) \eqs \{0\}$ if $\|u \|_{\Delta} < 1$.
Thus if the maximizer occurs at an interior point of $C_1$ then $m(b,v) = 0 \foral v \in \HzLap$ and the maximum value is 
0 which is not true.
Thus the maximizer occurs at a $b$ with $\| b \|_{\Delta} \eqs 1$ and then
\[ \pal I_1 (b) \eqs \{ \mu \, D a(b, . ) : \mu \leq 0 \} \]
as in the proof of proposition 9.6.1. of \cite{ABM}. Thus the maximizers  $b$ satisfy
\beq \label{e5.7}
\mu a(b,v)  \eqs m(b,v) \foral v \in \HzLap \ \mbox{and some} \quad \mu \geq 0. \eeq
Put $v = b$ to see that $\mu$ will be this maximum value $\beta_1$ and then \eqref{e5.7} shows that a maximizing 
$b_1$ satisfies  \eqref{e5.1} with $q_1 = {\beta_1}^{-1}$.   
Moreover $q_1$ will be the smallest eigenvalue of \eqref{e5.1} and the inequality \eqref{e5.5} holds by scaling
the constraint.
\epf

Given this first DBS eigenvalue and eigenfunction, a family of successive eigenvalues and eigenfunctions 
is now  constructed sequentially. 
Suppose that the set $\{q_1, \ldots , q_{k-1}\}$  of (k-1) smallest eigenvalues of \eqref{e5.1} and a corresponding sequence
of $\Delta-$orthonormal eigenfunctions $\{b_1,\ldots, b_{k-1}\}$ has been found.
Let $V_k$ be the subspace spanned by this finite set of eigenfunctions.

Define $W_k \deqs \{ u \in \HzLap : a(u,b_j) = 0 \xfor 1 \leq j \leq k-1 \}, C_k \deqs C_1 \cap W_k$.
Consider the problem of maximizing 
$\Mcl(u)$ on $C_k$ and evaluating $\beta_k \deqs  \sup_{u \in C_k} \ \Mcl(u)$.
This problem has maximizers that provide the next smallest eigenvalue and  associated normalized 
eigenfunctions as described next.

\btm \label{T5.2}
Assume (B2) holds and the smallest (k-1) DBS eigenvalues $q_j$ are known with an associated family of 
$\Delta-$orthonormal eigenfunctions $b_j$. 
When $C_k$ as above,  there are functions $\pm b_k  \in C_k$ that maximize $\Mcl$ on  $C_k$.
 $b_k$ is a non-trivial solution of \eqref{e5.1} associated with the next smallest positive eigenvalue 
$q_k = 1/{\beta_k}$ of the DBS eigenproblem. Also  $a(b_k,b_j)  = m(b_k,b_j) = 0$ for all $ 1 \leq j \leq k-1$ and
$a(b_k,b_k) = 1$. 
\etm   \bpf
For each finite k, $C_k$ is non-empty, closed, convex and bounded in $\HzLap$. 
Hence the weakly continuous functional $\Mcl$ attains a finite maximum $\beta_k$ on $C_k$. 
$\beta_k > 0$ as there are infinitely many independent functions in $\HzLap$ with $\Dnu u \neq 0$ on $\bdy$.

Just as in the previous proof, the maximizers must obey $a(b_k, b_k) =1$ by homogeniety.
From lemma 4.1 in \cite{Au4} and the analysis of section 9.6 of \cite{ABM}, a maximizer of $\Mcl$ on $C_k$
satisfies the equation
\beq 
m(b, v) \eqs  a(w,v) \pls \mu a(b, v) \ \mbox{for some  $w$ in} \  V_k, \mu \geq 0 \ \mbox{and all} \ v \in \HzLap. \eeq
Put $v = b_j$ here, then $a(w,b_j)$ as $b_k  \in W_k$. 
Thus $w = 0$ and $b_k$ is  a solution of $a(b,v) = q \, m(b,v)$ for all $v \in \HzLap$.
Thus $\pm b_k$ is a solution \eqref{e5.1} associated with the next smallest positive eigenvalue 
$q_k = 1/{\beta_k}$ of the DBS eigenproblem. 
By construction $a(b_k,b_j)  = m(b_k,b_j) = 0$ for all $ 1 \leq j \leq k-1$, so the theorem holds.
\epf

This theorem says that there one can find  a countable sequence of $\Delta-$orthonormal DBS eigenfunctions 
$\Bc \deqs  \{b_k : k \geq 1 \}$  with each eigenfunction maximizing $\Mcl$ on a set $C_k$ as above. 
Let $\BHarm$  be the  subspace of all biharmonic functions in $\HzLap$. 
It is closed in view of the definition  via \eqref{e5.3}. 
The following result is an analog of parts of theorem 3.18 in \cite{GGS}. 
See also Ferrero, Gazzola and Weth \cite{FGW}.

\btm \label {T5.3}
Assume that $\Om$ satisfies (B2) and $\Bc$ is a  sequence of DBS eigenfunctions constructed by the above algorithm. 
The corresponding DBS eigenvalues $  q_j $ each have finite multiplicity and increase to $\infty$.
 $\Bc$ is a  $\Delta-$orthonormal basis of the subspace $\BHarm$ of $\HzLap$. \etm
\bpf
When $\Bc$ is constructed as above, it converges weakly to zero in $\HzLap$ as it is $\Delta-$orthonormal.
Thus $\Mcl(b_k) = \beta_k$ converges to zero as it is weakly continuous - or $q_k$ increases to $\infty$.
Hence each eigenvalue has finite multiplicity.  

Let V be the closed subspace spanned by $\Bc$. It will be a subspace of $\BHarm$ since each $b_k \in \BHarm$. 
If $v \in \BHarm$ is $\Delta-$orthogonal to V and $\Mcl(v) > 0$, then $\vt \deqs v / {\sqrt{a(v,v)}} $ has $a(\vt,\vt) =1$
and $\Mcl(\vt) > \beta_K$ for some large K. 
This contradicts the definition of $\beta_K$, so we must have $\Mcl(v) = 0$ for all  that are  $\Delta-$orthogonal to V.
The uniqueness of solutions of the Dirichlet biharmonic problem on regions obeying (B2) then yields that such a 
$v$ must be zero, so $\Bc$ is a maximal   $\Delta-$orthonormal set in $\BHarm$ as claimed. 
\epf

This theorem implies that a   biharmonic function $b$  has the spectral, or eigenfunction, representation
\beq \label{e5.12}
b(x) \eqs \sum_{j=1}^\infty \ \ang{b,b_j}_{\Delta}  \ b_j(x) \xon \Om \eeq
that converges in the $\Delta-$norm to $b$ from the basic representation theorem for vectors in a Hilbert space.

Define $\HzzLap$ to be the class of all functions in $\HzLap$ that also have $\Dnu u = 0$  on  $\bdy$. 
Since $\Dnu$ is a continuous linear map, $\HzzLap$ is a closed subspace of $\HzLap$.
The following results lead to  an   orthogonal decomposition that  is analogous to that  described in theorem 3.19 of
\cite{GGS} as well as to the decomposition \eqref{e3.7} given above. Namely 
\beq \label{e5.14}
\HzLap \eqs \HzzLap \oplus_{\Delta} \BHarm \eeq
 where $\oplus_{\Delta} $ indicates the orthogonal complement $\wrt$ the $\Delta-$inner product.

To see this, assume $u \in \HzLap$ and $\Dnu u = \eta $ on $\bdy$. 
Define $\Keta$ to be the affine subspace of $\HzLap$ of functions with $\Dnu u = \eta $ on $\bdy$.

\begin{lem} \label{L5.4}
When (B2) holds, $\eta$ as above, there is a unique  function $\bt \in \BHarm$ that minimizes $\| u \|_{\Delta}$ on 
$\Keta$.
It is a solution of 
\beq \label{e5.15}
\IOm \Delta b \ \Lapv \ dx \eqs 0 \foral v \in \HzzLap. \eeq.  \end{lem}
\vspace{-2em}
\bpf
Let $\Ac(u) := a(u,u)$ as in \eqref{e5.4} and consider the problem of minimizing $\Ac$ on $\Keta$.
$\Ac$ is strictly convex, coercive and weakly l.s.c. on $\HzLap$ as it is a norm and thus there is a unique minimizer $\bt$
of $\Ac$ on $\Keta$.
The extremality condition satisfied by $\bt$ is that $a(\bt,v) = 0$ for all $v \in \HzzLap$ so \eqref{e5.15} holds. 
\epf

Suppose that $P_B: \HzLap \to \BHarm$ is the linear map defined by $P_B u = \bt$ when $u \in \HzLap$ has
$\Dnu u = \eta$. 
Define $P_{00} u \deqs u \mns P_B u$, then $ P_{00} u \in \HzzLap$.
These  are complementary projections of $\HzLap$ to itself and \eqref{e5.14} is an orthogonal decomposition since 
\eqref{e5.15}  holds.

Section 3.3 of \cite{GGS} describes the DBS spectrum and eigenfunctions for the unit ball in $\RN$  explicitly - 
and the formula for arbitrary balls may then be found using scaling arguments. 
It would be of great interest to have further information about these eigenfunctions and eigenvalues for simple two and three
dimensional regions. 
There has been some computation of such eigenvalues starting with the work of Sigillito and Kuttler described in \cite{KS}.

\vspace{1em}
\section{Orthonormal Bases and Reproducing Kernels for $\Btom$.  } \label{s6}

The Bergman space $\Btom$ is the space of  weakly harmonic functions in $\Lt$ - that is, those that satisfy \eqref{e3.2}. 
Chapter 8 of Axler, Bourdon and Ramey \cite{ABR} provides an introduction to Bergman spaces  and early results regarding
 these spaces are  described in Bergman \cite{Ber} and Bergman and Schiffer  \cite{BS}.
In this section  the orthogonal projection of $\Lt$ onto $\Btom$  and an explicit $L^2-$orthonormal  basis of $\Btom$
will be described. 
This leads to an explicit formula for the Reproducing Kernel of $\Btom$ in terms of the sequence of DBS eigenfunctions generated 
as  in the preceding section.

The orthogonal projection of $\Lt$ onto $\Btom$ may be found by looking at a variational principle for the orthogonal complement.
Given $f \in \Lt$, consider the minimum norm variational problem of minimizing the functional $\Dc :  \HzzLap \to \Rb$ 
defined by
\beq \label{e6.1}
\Dc (\psi) \deqs \IOm \ | \, \Lap \psi  \mns f \, |^2 \ \ dx .  \eeq

\begin{lem} \label{L6.1}
 \quad Suppose that (B2) holds and $f \in \Lt$, then there is a unique minimizer $\tilde{\psi}$ of $\Dc$
 on $\HzzLap. \ \psit$ satisfies 
\beq \label{e6.3}
\IOm (\Lap \psit \mns f) \ \Lap \chi \ dx \eqs 0 \foral   \chi \in \HzzLap. \eeq
 \end{lem} 
\bpf 
As in the proof of lemma \ref{L5.4}, $\Dc$ is continuous, strictly convex and coercive on $\HzzLap$. 
Hence there is a unique minimizer of $\Dc$ on $\HzzLap$.
$\Dc$ is G-differentiable and the standard extremality condition implies \eqref{e6.3}.
\epf

The system \eqref{e6.3} is the weak version of the boundary value problem  
\[ \Delta^2 \, \psi \eqs \Delta \, f \xon \Om \xwith \ \psi \eqs \Dnu \psi \eqs 0 \xon \bdy. \]
The solution $\psit$ will be called the {\it biharmonic  potential} of $f$.
Define $P_W : \Lt \to \Lt$ by $P_W f \deqs \Lap \psit$.
It is straightforward to  verify that $P_W$ is a projection onto a subspace $W := \Delta(\HzzLap) \subset \Lt$. 
The  range of this $P_W$  is closed from corollary 3.3 of \cite{Au}.  

Define $\PcH \deqs I \mns P_W$ then $\PcH$ is also be a projection on $\Lt$ with closed range that will be called the
{\it Bergman harmonic projection}. 
The range of $\PcH$ is the class of all functions $v \in \Lt$ that satisfy
\beq \label{e6.4}
\IOm  v \, \Delta \chi \ dx \eqs 0 \foral \chi \in \HzzLap. \eeq
Thus $v $ is harmonic on $\Om$ as \eqref{e3.2} holds we have an $L^2-$orthogonal decomposition
\beq \label{e6.5}
\Lt \eqs \Btom \oplus W.  \eeq
This decomposition is a version  of a result  attributed to Khavin described in  lemma 4.2 of Shapiro \cite{Sha}.

The Bergman harmonic projection $\PcH$ is not the same as  the  harmonic projection of lemma \ref{L3.1}.
$P_H f$  the closest harmonic function to $f$ in the $\pal-$norm on $\Hone$ while $\PcH f$ is the closest harmonic function
 in  the $L^2-$ norm on $\Om$. 
In particular the standard  projection has $P_H f = 0$ for all $f \in \HzLap$ while $\PcH f$ may be non-zero for functions in
$\HzLap$.

 Let $\Bc \eqs \{b_j : j \geq 1 \}$ be a maximal $\Delta-$orthonormal sequence of DBS eigenfunctions 
constructed using the algorithm of  section \ref{s5}.
Define $h_j \deqs \Delta b_j$ for each $j \geq 1$ and   $\BH \deqs \{h_j : j \geq 1 \} \subset \Lt$.
The $h_j$  are harmonic from \eqref{e5.1} as 
\[ \IOm \ h_j \ \Delta v \ dx \eqs \IOm \ \Delta b_j \ \Delta v \ dx \eqs 0 \foral v \in C_c^2(\Om) \]
so \eqref{e3.2} holds as $C_c^2(\Om) \subset \HzzLap$.  
Since these $b_j$ are $\Delta-$orthonormal, the $h_j$ are $L^2-$orthonormal.

For $M \geq 1$, consider the functions $H_M : \Omt \to \R$ defined by 
\beq \label{e6.6}
H_M(x,y) \deqs \sumM \ h_j(x) \, h_j(y) \eqs \sumM  \Delta b_j(x) \,  \Delta b_j(y). \eeq
Since each $h_j$ is $C^{\infty}$ on $\Om$ from Weyl's lemma for harmonic functions, $H_M$ is also $C^{\infty}$ 
and the integral operator $\HcM : \Lt \to \Lt$   defined by 
 \beq \label{e6.7}
\HcM f(x) \deqs \IOm \, H_M(x,y) \, f(y) \, dy \eeq
is a finite rank projection of $\Lt$ into $\Btom$.

\btm \label{T6.1}
Assume (B2) holds and $\BH$ is constructed as above, then $\BH$ is a maximal orthonormal set in $\Btom$. 
The orthogonal projection $\PcH$ of $\Lt$ onto $\Btom$ has the representation
\beq \label{e6.8}
\PcH f (x) \eqs \lim_{M \to \infty} \ \HcM f(x) \eqs \sumo \ \ang{f,h_j} \, h_j(x) \foral f \in \Lt. \eeq  \etm
\bpf
The above comments show that $\BH$ is an orthonormal subset of $\Btom$. 
It remains to show that it is maximal.
Suppose not, then there is a $k \in \Btom$ with $\ang{k, h_j} = 0 $ for all $j \geq 1$.
Let $\ut $ be the unique solution in $\Hzone$ of the equation 
\[ \IOm  \gradu \cdot \gradv \ dx \eqs \IOm k \, v \, dx \foral v \in \Hzone. \]
This $\ut$ exists, is unique and  $- \Delta \ut = k \in \Lt$.
Thus $\ut \in \HzLap$ and $\IOm \Delta \ut \, \Delta b_j \, dx = 0$ for all $j \geq 1$.
So $k = \Delta \ut$ is $\Delta-$orthogonal to $\Btom$ from theorem \ref{T5.3}.
This implies $k = 0$ so $\BH$ is an orthonormal basis of $\Btom$.
Given that $\BH$ is an orthonormal basis of $\Btom$, the last sentence follows from the Riesz-Fisher theorm.
\epf

\begin{cor} \label{C6.3}
 \quad Suppose that (B2) holds, then $\Btom$ is a Reproducing Kernel Hilbert (RKH) space with  reproducing kernel
\beq \label{e6.9}
\ROm(x,y) \deqs \sum_{j=1}^{\infty} \ \Delta b_j(x) \ \Delta b_j(y) \xfor (x,y) \in \Omt. \eeq 
 \end{cor} 
 \vspace{-1em}
\bpf 
The fact that $\Btom$ is a RKH space is classical; see theorem 8.4 of \cite{ABR}.
From  theorem \ref{T6.1}, when  $k \in \Btom$ then
\beq \label{e6.10}
k(x)  \eqs \sum_{j=1}^{\infty} \ \kjh  h_j(x) \eqs \sum_{j=1}^{\infty} \  \ang{k, \Delta b_j} \ \Delta b_j(x) \xon \Om \eeq
This  series converges in the $L^2-$norm so  $k(x) \eqs \ang{\ROm(x,.), k}$ when  $x \in \Om$ and $\ROm$ is defined by 
\eqref{e6.9}.

 \epf

It appears that this expression is an eigenfunction expansion of the very different  reproducing kernel found by 
J.L. Lions in \cite{L1}. 
Note in particular that the reproducing kernel  $\ROm$ acts as a {\it Delta function} on the subspace $\Btom$ of $\Lt$
in that 
\beq \label{eDel}
k(x) \eqs \IOm \ \ROm(x,y) \ k(y) \ dy \foral  x \in \Om \xand  k \in \Btom. \eeq

When $b_j$ is a DBS eigenfunction, define $w_j := \sqrt{q_j |\bdy|}  \ \Dnu b_j$. 
From theorem \ref{T3.1}  $w_j \in \Ltby$. 
The boundary conditions on the DBS eigenfunctions imply that
\beq \label{e6.11}
\gamma(h_j) \eqs q_j  \ \Dnu b_j \eqs \sqrt{\frac{q_j}{|\bdy|} } \, w_j \foreach j \geq 1. \eeq

Let $\Wc \deqs  \{ w_j : j \geq 1 \}$, then the following result says $\Wc$ is an orthonormal  basis of $\Lttby$.

\btm \label{T6.2}
Assume (B2) holds,  then $\Wc$ is a maximal orthonormal set in $\Lttby$. \etm
\bpf
From the definition \eqref{e5.1}, the DBS eigenfunctions satisfy
\[ \ang{\Delta b_j, \Delta b_k} \eqs q_j \ |\bdy| \  \ang{\Dnu b_j, \Dnu b_k}_{\bdy} \]
for all j,k. 
Since the $\Delta b_j$ are orthonormal in $\Btom$,  the $w_j$ are orthonormal in $\Ltby$. 

Suppose these functions are not maximal. Then there is a non-zero  $\eta \in \Lttby$ such that 
$\ang{\eta, w_j}_{\bdy} = 0$ for all $j \geq 1$.  
Let $\htil  = E \eta \in \Hharm$ be the harmonic extension of this boundary data defined as in section  \ref{s4}.
This function is in $\Btom$ so $\htil \neq 0$ and it has a representation of the form \eqref{e6.8}. 
In this case at least one $\ang{\htil, h_j} \neq 0 $.
In view of \eqref{e6.11} this contradicts our assumption, so $\Wc$ is an orthonormal basis of $\Lttby$.
\epf

In view of this result and the fact that $\BH$ is an orthonormal basis of $\Btom$, we shall regard the 
{\it harmonic trace operator} $\gamH$ to be the linear transformation  on $\Btom$ defined  by 
\beq \label{e6.12}
\gamH(k)   \eqs {|\bdy|}^{-1/2} \ \sumo \sqrt{q_j} \, \hat{k}_j w_j \ \mbox{   when} \quad  k = \sumo \hat{k}_j h_j  \in \Btom. \eeq
From \eqref{e6.11} this is an unbounded map into $\Ltby$.
In Auchmuty \cite{Au3}, it was shown that the harmonic  trace operator  defines an isometric isomorphism between $\Lt$ and a space 
that was denoted   $\Hmhby$ and characterized as the dual space of $\Hhby \ \wrt$ the usual inner product on 
$\Lttby$.

\vspace{1em}
\section{SVD Representation of the Poisson Kernel  } \label{s7}

The preceding analysis of the Dirichlet Biharmonic Steklov eigenfunctions showed that  they generate orthonormal bases
of the spaces $\BHarm, \Btom$ and $\Lttby$.
A closer investigation shows that they also provide a {\it singular value decomposition} (SVD) of the usual Poisson 
integral operator  for solving the Dirichlet problem for harmonic functions. 

Given a function $g \in \Ltby$, consider  the harmonic extension problem  of a  finding  $\ut := E_H g \in \Btom$ 
satisfying  $\gamH(\ut) = g$. 
Then $E_H: \Ltby \to \Btom$ is a right  inverse of $\gamH$ and \eqref{e6.11} shows that, when 
$ g = \sum_{j=1}^{\infty} \gjh \, w_j \in \Lttby$, then  
\beq \label{e7.5}
E_H g(x) \eqs  \sqrt{|\bdy|} \ \sum_{j=1}^{\infty} {\frac{\gjh}{\sqrt{q_j}}} \   h_j(x)   \eeq

\btm \label{T7.1}
Assume that (B2) holds and  $E_H: \Ltby \to \Btom$ is defined by \eqref{e7.5}.
Then $E_H$ is an injective compact linear transformation with $\|  E_H  \| = \frac{1}{\sqrt{q_1}} $.   
\etm
\bpf
From \eqref{e7.5} and  orthonormality, Parseval's equality  yields
\[ \| \, E_H g \, \|^2 \eqs |\bdy|  \sumo   \ \frac{1}{q_j} \ \gjh^2 \leqs  \frac{|\bdy|}{q_1} \ \| g \|^2_{\bdy}. \]
so $E_H$ is continuous with  norm  as in the theorem.
$E_H$ is obviously injective. It is compact as the $q_j$ increase to $\infty$. 
\epf

This formula for the norm of $E_H$  is well-known when the boundary obeys stronger boundary regularity conditions. 
It was  described in Fichera \cite{Fic} and is  equation 3.1 in \cite{Sig}. 
Note that \eqref{e7.5} is essentially a  SVD of this harmonic extension operator as it maps one orthonormal basis to another.
The singular values of the operator are simple functions of   the  eigenvalues of the Dirichlet Biharmonic Steklov problem.

Classically this operator has usually been described in terms of the {\it Poisson kernel} - see section 2.2 of Evans \cite{Ev}
or most  other PDE texts. The solution of this boundary value problem is  described in terms of a function   
$P: \Om \times \bdy \to [0,\infty]$ such that  
\beq \label{e7.7}
E_H g(x) \deqs \Iby \, P(x,z) \, g(z) \dsg. \eeq

A comparison of \eqref{e7.7} and \eqref{e7.5} leads to the following singular value decomposition of the 
Poisson kernel as a function on $\Om \times \bdy$.

\btm \label{T7.2}
Assume that (B2) holds, then  the Poisson kernel $P(x,z)$ has the  singular value   representation
\beq \label{e7.8}
P(x,z) \eqs |\bdy|^{-1/2} \  \sum_{j=1}^{\infty}  \ \frac{1 }{\sqrt{q_j}} \ h_j(x) \,  w_j(z) \eqs  \sum_{j=1}^{\infty} \, \frac{1}{q_j} \  \Delta b_j(x) \, \Delta b_j(z) \eeq
for $(x,z) \in \Om \times \bdy$. 
\etm
\vspace{-1em}  \bpf
The formulae in \eqref{e7.8} hold by comparing \eqref{e7.7} and \eqref{e7.5} and using the definitions and properties of the various functions. 
 \epf

The singular value decomposition of theorem \ref{T7.1} leads to explicit formulae for the best rank M 
approximations of the Poisson operator. 
For finite M, define $P_M$ and $E_M$ by
\beq \label{e7.9}
P_M(x,z) \deqs \sumM \ \frac{h_j(x)}{\sqrt{|\bdy| \, q_j }} \ w_j(z) \xand E_M g(x) \deqs \Iby P_M(x,z) g(z) \dsg(z). 
\eeq

Then for each $z \in \bdy, \ P_M(.,z)$ is a harmonic function on $\Om$ and for each $x \in \Om, \ P_M(x,.)$
is an $L^2-$function on $\bdy$ with 
\[ \| P_M(x,.)\|_{\bdy}^2 \eqs \sumM \frac{|h_j(x)|^2}{|\bdy| \, q_j} \xand \IOm \Iby |P_M(x,z)|^2 \dsg \, dx \eqs \sumM \frac{1}{q_j}. \]

These formulae lead to the following approximation result for the Poisson  operator. 
\btm \label{T7.3}
Assume that (B2) holds, $E_H$ is the harmonic extension  operator and $P_M, E_M$ are defined 
by \eqref{e7.7}, then
\beq \label{e7.11}
\| \, E_H g - E_M g \, \| \leqs \sqrt{\frac{|\bdy|}{q_{M+1}}} \ \| \, g - g_M \, \|_{\bdy} \foral g \in \Ltby \eeq
Here $q_{M+1}$ is the $(M+1)-$st DBS  eigenvalue and $g_M := \sumM \gjh w_j$.
\etm
\bpf
From \eqref{e7.7} and \eqref{e7.9}, one sees that
\[ E_H g(x) \mns  E_M g(x) \eqs |\bdy|^{1/2}  \, \sum_{j=M+1}^{\infty} \ {\frac{\gjh}{\sqrt{q_j}}} \   h_j(x), \quad \mbox{so} \]
\[ \|  E_H g -  E_M g  \|^2 \eqs  |\bdy|  \sum_{j=M+1}^{\infty} \  \frac{1}{q_j} \ \gjh^2 \leqs  \frac{|\bdy|}{q_{M+1}} \ \| g- g_M \|^2_{\bdy}  \]
so \eqref{e7.11} holds as claimed.
\epf

   \vspace{1em}
 \section{Normal Derivatives  of Laplacian Eigenfunctions} \label{s8}

Theorem \ref{T3.1} provided an estimate of the normal derivative of functions in $\HzLap$. 
Here a result about  the normal derivatives of  eigenfunctions of the Dirichlet Laplacian will be proved. 
This quantifies part of theorem 1.1 of Hassell and Tao \cite{HT}  that answered a question of Ozawa \cite{Oza}.

A non-zero function $e \in \HzLap$ is said to be a Dirichlet eigenfunction of the Laplacian on $\Om$
corresponding to an eigenvalue $\lambda$ provided 
\beq \label{e8.1}
\IOm [ \ \nabla e \cdot \gradv \mns \lambda \ e \, v \ ] \ dx \eqs 0 \foral v \in \Hzone. \eeq

The eigenfunction is normalized if $\| \, e \, \| = 1$. 
Note that theorem \ref{T3.1} already provides a generic  upper bound for the constant in the inequality 1.1 of \cite{HT}.
Here an explicit representation  for the normal derivatives of Dirichlet eigenfunctions will be derived
that shows this constant may be bounded in terms of the first DBS eigenvalue $q_1$.  

\btm \label{T8.1}
Suppose (B2) holds and $e \in \HzLap$ is a normalized Dirichlet eigenfunction of the Laplacian on $\Om$ with 
eigenvalue $\lambda$.   Then
\beq \label{e8.3}
\Iby |\, \Dnu e \, |^2 \dsg  \leqs  \ \frac{\| \PcH e \|^2}{q_1} \ \lambda^2  \leqs \frac{1}{q_1} \ \lambda^2 \eeq 
with  $\PcH$  the Bergman harmonic projection. \etm
\bpf
When $\psi$ is the biharmonic potential of $e$ then, from \eqref{e6.5}, one has $e = \Delta \psi + e_H$ with
$\psi \in \HzzLap$ and $e_H = \PcH e \in \Btom$. 
From theorem \ref{T6.1}, $e_H = \sumo \, \ejh \Delta b_j$, so the eigenvalue equation yields that
\[    \Delta \ ( \, \lambda^{-1} e \pls \psi \pls \sumo \, \ejh  b_j \, ) \eqs 0 \xon \Om. \]
The function here is also zero on $\bdy$, so
\beq \label{e8.5}
e(x) \ \equiv \ - \lambda \, [ \, \psi \pls \sumo \, \ejh  b_j \, ] \xon \Om. \eeq
Thus $\Dnu e \eqs - \lambda \,  \sumo \, \ejh \,  \Dnu b_j$ on $\bdy$ as  $\psi \in \HzzLap$.
The definition of $w_j$ implies that
\beq \label{e8.7}
	\Dnu e(x) \eqs - \lambda \,  \sumo \, \frac{\ejh}{\sqrt{|\bdy| \, q_j}} \,  w_j(x) \xon \bdy. \eeq
Since $q_1$ is the least DBS eigenvalue, \eqref{e8.3} now follows from Parseval's equality as
$\|  \PcH e \|^2 = \sumo \ejh^2$ and $\| \PcH \| = 1$.
\epf

In particular this result  shows that the constant C in the Hassell-Tao inequality  for nice bounded regions in
$\RN$ has $C \leq 1/\sqrt{q_1}$.   

These eigenfunctions illustrate a difference between the Steklov harmonic projection which has $P_H e = 0$
 for any Dirichlet   eigenfunction, and  the Bergman harmonic projection $\PcH$ which must have $\PcH e \neq 0$.  

\vspace{1em}
%

\end{document}